\documentstyle[amsfonts,amssymb]{article}

\title {On generalized Kostka polynomials \\ and quantum Verlinde rule}
\author{B.~Feigin, S.~Loktev\thanks{Partially supported by CRDF grant RM1--265}}
\date{}
\newcommand{\name}[1] {\mathop{\rm #1}\nolimits}
\newcommand{\itname}[1] {\mathop{\it #1}\nolimits}

\newcounter {nst}
\newenvironment {statement}[1] {\refstepcounter{nst}\par\medskip%
\noindent{\bf #1 \, \arabic{section}.\arabic{nst}. }\it} {\par}
\newenvironment {Statement}[1] {\par\medskip\noindent{\bf
#1. }\it}{\par\medskip}
\newcommand{\sbros}{\setcounter{nst}{0}}

\newcommand \qed {\,\rule[-.23ex]{1.6ex}{1.6ex}}
\newcommand{\proofName}{Proof.}
\newenvironment {proof} {\par\medskip\noindent{\bf \proofName\,}}
{\qed\par\medskip}

\newcommand{\gref}[2]{\ref{#1}.\ref{#2}}

\newcounter{tmp}
\newcommand{\rf}{{\par\medskip\setcounter{tmp}{1} ${\rm (\roman{tmp})}$\
}}
\newcommand{\rr}{{\par\medskip\refstepcounter{tmp} ${\rm (\roman{tmp})}$\
}}

\newcommand{\gtc}[1]{{\frak #1}}
\newcommand{\bbb}[1]{{\Bbb #1}}

\newcommand{\toki}[1]{{\frak #1}^{\Bbb C}}
\newcommand{\tg}{{\widehat{\gtc{g}}}}
\newcommand{\chal}{{\check{\alpha}}}
\newcommand{\htg}{{\gtc{h}_\tg}}
\newcommand{\ver}[1]{{{\cal{V}}_{#1}(\gtc{g})}}
\newcommand{\subw}{{{w \in W_{aff}}\atop{w(\lambda,k) \ \mbox{\scriptsize
is dominant}}}}
\newcommand{\conf}{{\bbb{C}^n \setminus \Delta}}
\newcommand{\ch}{\mbox{ch}}
\newcommand{\gr}{\mbox{gr}^\bullet}
\newcommand{\nmod}{{\widetilde{{\cal{V}}}(\gtc{g})}}

\newcommand{\tsl}{{\widehat{sl_2}}}

\newcommand{\np}{{{\gtc{n}}}}
\def \A {{\frak A}}
\def \Complex {{\Bbb C}}
\def \sltwo {\itname{sl}_2}

\begin{document}

\maketitle

\begin{abstract}
Here we propose a way to construct generalized Kostka polynomials. Namely,
we construct an equivariant filtration on tensor products of irreducible
representations. Further, we discuss properties of the filtration and
the adjoint graded space. Finally, we apply the construction to 
computation of coinvariants of current algebras.
\end{abstract}

\section*{Introduction}
\def \slhat {\widehat{\sltwo}}

In this paper we consider the following question. Let $\pi$ be a
representation from a certain minimal conformal field theory. For example
$\pi$
may be an irreducible representation of the Virasoro algebra from a
certain $(p,q)$--theory (see \cite{BPZ}). Choose the standard basis
$\{L_j\}$ of the Virasoro
algebra and suppose that the vacuum vector in $\pi$ is annihilated by
$L_j,\, j < 0$. Let $A_n$ be the subalgebra formed by $L_n, L_{n+1},
\dots$.
Consider the space of coinvariants $\pi/A_n\pi$. In \cite{BFM} it is
proved that
$\dim (\pi/A_n\pi) < \infty$. The space $\pi/A_n\pi$ has the natural
grading,
$\pi/A_n\pi = \bigoplus (\pi/A_n\pi)^j$ such that $\deg L_i = i$ and the
vacuum vector
has grading
zero. So define the truncated character $\name{ch}_q(\pi,n)$ as the
polynomial $\sum q^j \dim \, (\pi/A_n\pi)^j$. The problem is how to find
$\name{ch}_q(\pi,n)$.

First of all, let us try to evaluate the dimension of $\pi/A_n\pi$, that
is, $\name{ch}_1(\pi,n)$. Let $n$ be odd, $n = 2k-1$. The algebra
$\{L_{-1}, L_0, L_1, \dots\}$ can be considered as the algebra of
polynomial vector fields
on the line $f(t)\, \partial/\partial t$. Then
the subalgebra $A_{2k-1}$ consists of the fields $t^{2k}f(t)\,
\partial/\partial t$. Consider the family of Lie algebras
$$A_{2k-1}(p_1,
\dots, p_k) = (t-p_1)^2 \dots (t-p_k)^2 \Complex[t]\, \partial/\partial
t,$$
$p_1, \dots, p_k$ --- points on the line. It is easy to prove that 
$$\dim
\pi/A_{2k-1}\pi \ge \dim \pi/A_{2k-1}(t_1, \dots, t_k)\pi$$ 
because
$A_{2k-1}(t_1, \dots, t_k)$ is a filtered deformation of $A_{2k-1} =
A_{2k-1}(0,0, \dots, 0)$. It is natural to conjecture that actually
%*
\begin{displaymath}
\dim \pi/A_{2k-1}\pi = \dim \pi/A_{2k-1}(t_1, \dots, t_k)\pi.
\end{displaymath}
%*
Still it is an open problem. In \cite{FFr} it was proved for
$(2,2s+1)$--minimal theories, and in \cite{FFu} for $k = 2$. Suppose
that all the
points $t_1, \dots, t_k$ are distinct. In this case the dimension of
$\pi/A_{2k-1}(t_1, \dots, t_k)\pi$ is given by the Verlinde formula
(see \cite{BFM}).
More precisely, let ${\cal{V}}_N$ be the Verlinde algebra, $\{ \pi_1,
\dots,
\pi_N \}$
be
the standard basis in ${\cal{V}}_N$ corresponding to the set of
``integrable''
representations. Let us write down in ${\cal{V}}_N$
{\def \theequation {*}
%*
\begin{equation}
(\pi_1 + \dots + \pi_N)^k = \sum_{j=1}^N c_j(k) \cdot \pi_j.
\end{equation}
%*
}
Then for distinct $\{t_j\}$ we have
%*
\begin{displaymath}
c_j(k) = \dim \pi_j/A_{2k-1}(t_1, \dots, t_k)\pi_j.
\end{displaymath}
%*

So in the case when the conjecture is true our problem is to find a
``$q$--version'' of the formula $(*)$. But, as we mentioned, the
conjecture is known only in $(2, 2s+1)$--case. In general case we can do
the following.
The representation $\pi$ has the natural grading: $\pi = \pi^0 \oplus
\pi^1
\oplus \dots$. Let $\pi^{\le j} = \pi^0 \oplus \pi^1 \oplus \dots \pi^j$.
Consider the composition $\pi^{\le j} \hookrightarrow \pi \to
\pi/A_{2k-1}(t_1,
\dots, t_k)\pi$ and let $S_j$ be the image of $\pi^{\le j}$. Therefore,
we have
a sequence of spaces $S_0 \subset S_1 \subset \dots$. Let us write down
the ``character'' of the space $\pi/A_{2k-1}(t_1, \dots, t_k)\pi$ as $\sum
\dim S_j/S_{j-1} \cdot q^j$. If the main conjecture is true, then this
character coincides with the character of $\pi/A_{2k-1}\pi$.

Now let us rewrite $(*)$ in the ``alternating sum'' form. 
Let $U$ be a ``small''
quantum group $U_q(sl_2)$ in the root of unity of order $l$. $U$ is
generated as an
algebra by $E, F, K, h$ with usual relations, in particular $E^l = F^l =
0$. Let $Q$ be a representation of $U$ where $h$ acts
by the diagonal way. Denote by $Q(\lambda)$ the subspace in $Q$ where $h$
acts by the scalar $\lambda$. The ``toy Felder complex'' $F(Q)$ is by
definition
%*
\begin{displaymath}
\dots \stackrel{E}{\gets} Q(-2l) \stackrel{E^{l-1}}{\gets} Q(-2)
\stackrel{E}{\gets} Q(0) \stackrel{E^{l-1}}{\gets} Q(-2+2l)
\stackrel{E}{\gets} Q(2l) \gets \dots
\end{displaymath}
%*

Let $\{\nu_0, \dots, \nu_{l-2}\}$ be the set of irreducible Weyl
modules of $U$, $\dim \nu_\alpha = \alpha+1$, let $\nu_{l-1}$ be the
Steinberg
module.

\begin{Statement}{Proposition}
Let $\nu_{\alpha_1}, \dots, \nu_{\alpha_r}$ be a set of Weyl modules, all
$\alpha_j < l-1$. The homology of the complex $F(\nu_{\alpha_1} \otimes
\dots \otimes \nu_{\alpha_r})$ are concentrated in the term $\nu_{\alpha_1}
\otimes \dots \otimes \nu_{\alpha_r}(0)$. The dimension of these homologies
can be calculated in the following way. Representations $\nu_0, \dots,
\nu_{l-2}$ form a basis in a certain Verlinde algebra. Write down
$\nu_{\alpha_1} \cdot \dots \cdot \nu_{\alpha_r} = \sum c_i \nu_i$. The
dimension of the homologies is equal to the coefficient $c_0$.
\end{Statement}

Let us apply this to the simplest $(2,5)$--theory. In this case there are
two representations in the minimal theory. Namely they are the vacuum
representation
$\epsilon_0$, and the representation $\epsilon_1$ with the highest weight
$-1/5$. The Verlinde algebra for $(2,5)$ is generated by
$\epsilon_0$ and $\epsilon_1$ with relations $\epsilon_0\cdot \epsilon_i
= \epsilon_i$, $i=0,1$ and
$\epsilon_1^2 = \epsilon_0 + \epsilon_1$.
The quantum group related to the $(2,5)$--theory is $U_q(\sltwo),\, q^5 =
1$.
In
terms of $U_q(\sltwo)$ the Verlinde algebra can be  described in the
following
way. $U_q(\sltwo)$ has four Weyl modules $\nu_0, \nu_1, \nu_2, \nu_3$
($\nu_4$ is Steinberg), and we have to impose an additional relation:
$\nu_0 = \nu_4$. After that $\nu_0$ corresponds to $\epsilon_0$ and $\nu_1$
to $\epsilon_1$. Note that $\epsilon_1^2 = \epsilon_0 + \epsilon_1$,
therefore $(\epsilon_0 + \epsilon_1)^{k} = \epsilon_1^{2k}$. Summarizing,
we can conclude that
$\dim \epsilon_0 / A_{2k-1} \epsilon_0$ coincides with Euler
characteristics of the
complex $F(\nu_1\otimes \dots \otimes \nu_1)$, $\nu_1$ appears $2k$
times.

If we consider the $(2,2s+1)$--theory then the similar statement is true.
The quantum group now is $U_q(\sltwo)$, $q^{2s+1} = 1$. 
Denote by
$\epsilon_0$ the vacuum representation, and denote by $\nu_{s-1}$ the
$s$--dimensional representation of $U_q(\sltwo)$. Then
$\dim \epsilon_0 / A_{2k-1} \epsilon_0$ coincides with Euler
characteristics of the
complex $F(\nu_{s-1} \otimes \dots \otimes \nu_{s-1})$, $\nu_{s-1}$
appears $2k$ times. 

There exists a filtration on the toy Felder complex $F(\nu_s \otimes \dots
\otimes \nu_s) = F_0 \supset F_1 \supset \dots$ such that each $F_i$ is a
subcomplex, and the corresponding spectral sequence degenerates in the
first term. It is clear that the homology of the complex $\bigoplus
F_i/F_{i+1}$ is a graded space, and we claim that this graded space  is
isomorphic to $\epsilon_0 / A_{2k-1} \epsilon_0$.

Construction of the filtration on the tensor product is rather
complicated, and we will not discuss it in this paper. But actually the
filtration (but not the structure of complex) exists on the tensor
products of irreducible representations of $U_q(\sltwo)$ for arbitrary $q$,
even in the ``classical'' case $q = 1$. In this paper we consider the
classical situation but not only in the $\sltwo$--case. Note that from a
different point of view the problem of ``truncated'' character formula was
studied (and is actually being studied) in the articles of the Stony Brook
group (\cite{BMS}, \cite{SW1}, \cite{SW2}) and many others (\cite{FLOT},
\cite{KSz}, \cite{LT}, \cite{Sz}). Our work is
inspired by these remarkable papers.

In the first part of the paper we present the construction of the
filtration on the tensor product of finite--dimensional representations of
a Lie algebra (in principle arbitrary, but, actually, we need some
additional
structure on the representations). This construction is interesting by
itself, and we discuss in the second part its connection with the geometry
of the Schubert
variety, and with the structure of integrable representations of affine
Kac--Moody algebras. In the third part of the paper we discuss the problem
about the characters of coinvariants for the current algebras and present
a formula for $q$--Verlinde rule. 
The $\slhat$ situation of this formula is quite similar to the
$(2,2s+1)$--theories
for Virasoro.

\section{Fusion product of finite--dimensional
representations}\label{init}

\subsection{Filtration on the tensor product.}

Let $\gtc{a}$ be a Lie algebra. We call an $\gtc{a}$--module $\pi$ {\it
cyclic} if it contains a chosen vector $v$ such that $\pi = U(\gtc{a}) v$. 

Consider $n$ finite--dimensional cyclic $\gtc{a}$--modules
$\pi_1,\dots, \pi_n$. 
Our purpose is to define an $\gtc{a}$--equivariant  
filtration on their tensor product $\Pi = \pi_1 \otimes \dots \otimes
\pi_n$.

Of course, this definition requires some additional data. Namely, the
filtration depends on
a complex vector $ {\cal{Z}}=(z_1, \dots z_n)$. 
This dependence makes possible to distinguish isomorphic
submodules.

To construct the filtration, consider the algebra $\toki{a} =
\gtc{a}\otimes \Complex [t]$ of $\gtc{a}$--valued complex polynomials. 

Let $\pi$ be an $\gtc{a}$--module, $z \in \Complex$. Then we can define
an action of $\toki{a}$ on $\pi$. Namely, let $P \in \toki{a}$ be a
$\gtc{g}$--valued polynomial, $u \in \pi$. Then
$$P \cdot u = P(z) u.$$
This $\toki{a}$--module is called an {\it evaluation} representation.
Denote it by $\pi(z)$.

Now let $\pi_1, \dots, \pi_n$ be cyclic $\gtc{a}$--modules with chosen
cyclic vectors $v_1, \dots, v_n$, ${\cal{Z}} = (z_1, \dots z_n)$ be a
complex vector. Then we can define an action of $\toki{a}$ on
their tensor product $\Pi$. Namely, consider the $\toki{a}$--module

$$\Pi ({\cal{Z}}) = \pi_1(z_1) \otimes \dots \otimes \pi_n(z_n).$$

\begin{statement}{Proposition}\label{wdef1}
If $z_i$ are pairwise distinct numbers then the
$\toki{a}$--module
$\Pi({\cal{Z}})$ is generated by $v_1\otimes \dots \otimes v_n$.
\end{statement}

\begin{proof} 
Let $u_1 \otimes \dots \otimes u_n$ be a monomial in
$\Pi ({\cal{Z}})$,
$a \in \gtc{a}$. Applying the action of $a t^k \in \toki{a}$ to
this monomial we obtain that
$$a t^k(u_1\otimes \dots \otimes u_n) = \sum_i z_i^k \cdot 
u_1\otimes \dots \otimes a(u_i)\otimes \dots \otimes u_n.$$

As $z_i$ are pairwise distinct, the matrix $(z_i^j)$ is invertible, hence 
we obtain any monomial  $u_1\otimes \dots \otimes a(u_i)\otimes \dots
\otimes
u_n$ as a linear combination of $a t^j(u_1\otimes\dots\otimes u_n)$.

Now consider the vector $v_1\otimes \dots \otimes v_n$. Iterating the
procedure described above with different elements of $\gtc{a}$ we obtain
that
$$U(\gtc{a})v_1\otimes \dots \otimes U(\gtc{a})v_n \subset U(\toki{a})(v_1
\otimes \dots \otimes v_n).$$

As $v_i$ are cyclic vectors, this complete the proof.
\end{proof}

Denote the set of $n$ pairwise distinct complex numbers by $\conf$.
So if ${\cal{Z}} \in \conf$ then $\Pi({\cal{Z}})$ is a cyclic
$\toki{a}$--module.

Note that the Lie algebra $\toki{a}$ is graded by the degree of
polynomial. Then the algebra $U(\toki{a})$ is also graded. Thus for any
cyclic $\toki{a}$--module $M$ with cyclic vector $v$ we have the
induced filtration
$$F^i M = U^{\le i}(\toki{a}) v.$$
In particular, $F^0 M$ is the $\gtc{a}$--submodule generated by $v$.

\begin{statement}{Definition}
Let ${\cal{Z}} \in \conf$ be a set of
pairwise distinct complex numbers, $\pi_i$ be cyclic $\gtc{a}$--modules.
Then the {\rm filtered tensor product} ${\cal{F}}_{\cal{Z}}(\pi_1, \dots,
\pi_n)$ is the space $\Pi({\cal{Z}}) \cong \Pi$ with the induced
filtration.
\end{statement}

\begin{statement}{Proposition}
This is a well--defined increasing
$\gtc{a}$--equivariant filtration on the
tensor product $\Pi = \pi_1\otimes \dots \otimes \pi_n$.
\end{statement}

An $\gtc{a}$--module settled in the origin $\pi(0)$ provides us the
simplest example of a graded $\toki{a}$--module. And we can similarly
define the filtered tensor product for graded cyclic $\toki{a}$--modules.

First of all, if $\pi$ is a graded $\toki{a}$--module, $z$ is a complex
number, then we can define an $\toki{a}$--module $\pi(z)$, such that
$$P(t) \cdot u = P(t+z) u$$
for $P \in \toki{a}$, $u \in \pi(z)$. 

As above, let $\pi_1, \dots , \pi_n$ be cyclic graded
$\toki{a}$--modules with cyclic vectors $v_1, \dots, v_n$, ${\cal{Z}} =
(z_1, \dots z_n)$ be a complex vector. Then consider the
$\toki{a}$--module

$$\Pi ({\cal{Z}}) = \pi_1(z_1) \otimes \dots \otimes \pi_n(z_n).$$

\begin{statement}{Proposition}
If $z_i$ are pairwise distinct numbers then the
$\toki{a}$--module $\Pi ({\cal{Z}})$ is generated by the vector $v_1
\otimes
\dots \otimes v_n$.
\end{statement}

\begin{proof}
Let $u = u_1\otimes \dots \otimes u_n$ be a monomial in
$\Pi ({\cal{Z}})$. We have:
$$a t^k(u) = \sum_{i,r} \left[{k}\atop{r}\right] z_i^{k-r} \cdot u_1
\otimes \dots \otimes a t^r(u_i) \otimes \dots \otimes u_n.$$ 

As $\pi_i$ are finite--dimensional and graded, there exists a sufficiently
large $N$ such
that $\gtc{a}\otimes t^k$ acts trivially for $k>N$. Consider the matrix
of the coefficients 
$$\left( \left[{k}\atop{r}\right] z_i^{k-r}\right)_{k,(i,r)}, \quad
i=1\dots n, \ r=1\dots N, \ k=1\dots nN.$$

Similar to the Vandermonde matrix, this matrix is invertible under the
conditions of the proposition. Thus  we
can obtain any resulting monomial as a linear combination of $at^k(u)$.
The rest of proof is similar to the proof of
proposition~\gref{init}{wdef1}.
\end{proof}

So we can introduce a similar definition.

\begin{statement}{Definition}
Let ${\cal{Z}} \in \conf$ be a set of
pairwise distinct complex numbers, $\pi_i$ be cyclic
graded $\toki{a}$--modules.
Then the {\rm filtered tensor product} ${\cal{F}}_{\cal{Z}}(\pi_1, \dots,
\pi_n)$
is the space $\Pi({\cal{Z}}) \cong \Pi$ with the induced filtration.
\end{statement}

Note that along the way we define an action of $\toki{a}$ on the tensor
product of $\gtc{a}$--modules and a non--standard action of $\toki{a}$ on
the tensor product of graded $\toki{a}$--modules.

\subsection{Case of simple Lie algebra.}

Now let $\gtc{g}$ be a simple Lie algebra. Choose a Borel subalgebra
$\gtc{b} \in \gtc{g}$ to determine highest vectors.

Any finite--dimensional $\gtc{g}$--module is a direct sum of irreducible
modules with highest vectors. One can easily show that any
direct sum of pairwise non--isomorphic irreducible
modules is a cyclic $\gtc{g}$--module. 

If $\pi$ is an irreducible representation then choose its
highest vector
as a cyclic vector. Indeed, in this case we can choose any vector, but the
filtration calculated with the chosen highest vectors is the most subtle.
Namely,
one can prove that if $\pi_1, \dots, \pi_n$ are irreducible
representations and $v_i$ are their highest vectors then for any $u_i \in
\pi_i$ and ${\cal{Z}}$
$$ {\cal{F}}^i_{\cal{Z}} ( (\pi_1, v_1), \dots, (\pi_n, v_n)) \subset
{\cal{F}}^i_{\cal{Z}} ( (\pi_1, u_1), \dots, (\pi_n, u_n)).$$

If $\pi$ is a sum of different irreducible
modules then we can consider the sum of highest vectors of its
components.

Now consider isomorphism classes of filtered finite--dimensional
$\gtc{g}$--modules. It is a discrete kind of data. On the other hand, all 
objects and transformations introduced above are algebraic. Therefore we
have 

\begin{statement}{Proposition}
For given $\pi_1, \dots \pi_n$ there exists a
Zariski open subset $U\in \conf$,
such that for any ${\cal{Z}}_1, {\cal{Z}}_2 \in U$ we have
$$ {\cal F}_{{\cal{Z}}_1}(\pi_1, \dots, \pi_n) \cong
{\cal F}_{{\cal{Z}}_2}(\pi_1,
\dots,\pi_n)$$
as filtered $\gtc{g}$--modules.
\end{statement}

We are convinced that $U$ always equals $\conf$. In 
section~\ref{fus} we discuss
a more general conjecture.

Note that this isomorphism is not natural.
If we fix the space $\pi_1\otimes \dots\otimes \pi_n$ then the filtration
has a non-trivial algebraic dependence on ${\cal{Z}}\in \bbb{C}$. 
Nevertheless, we can consider a Hilbert polynomial of this filtration
for a generic point ${\cal{Z}} \in U \subset \conf$.
$$\ch_q (\pi_1, \dots, \pi_n) = \sum q^i \cdot \ch \ \mbox{gr}^i
{\cal{F}}_{\cal{Z}}(\pi_1,\dots, \pi_n)$$

For irreducible representations we obtain a version of generalized Kostka
polynomials. We think that they coincides with the polynomials 
discussed in \cite{KSz}, \cite{LT}, \cite{SW2}.
For example if $\gtc{g} = sl_2$ then calculations leads to quantum
supernomial coefficients (\cite{SW1}). Also we have a coincidence
for some representations of $sl_n$.
We present the way of these calculations in section~\ref{count}.

For sums of irreducible representations we obtain more elaborate but
also very useful polynomials. In section~\ref{appl} we discuss the
relation
between such polynomials and the filtered Verlinde rule.

At last notice that in the case of a simple Lie algebra $\gtc{g}$ we have
$$ \gr {\cal{F}}_{\cal{Z}} (\pi_1, \dots, \pi_n) \cong \pi_1 \otimes \dots
\otimes \pi_n$$
as $\gtc{g}$--modules. This property motivates the following definition of
the fusion product.

\subsection{Fusion product}\label{fus}

Let $\gtc{a}$ be a Lie algebra, $\pi_1, \dots, \pi_n$ be
finite--dimensional graded cyclic $\toki{a}$--modules, ${\cal{Z}} \in
\conf$. Consider the $\toki{a}$--module ${\cal{F}}_{\cal{Z}} (\pi_1,
\dots,
\pi_n)$. Clearly, we have
$$ \gtc{a} t^k \cdot {\cal{F}}^i_{\cal{Z}} (\pi_1, \dots, \pi_n) \subset  
{\cal{F}}^{i+k}_{\cal{Z}} (\pi_1, \dots, \pi_n).$$ 

Hence, we have the well--defined action of $\toki{a}$ on the adjoint
graded space such that
$$ \gtc{a} t^k \cdot \mbox{gr}^i {\cal{F}}_{\cal{Z}} (\pi_1, \dots, \pi_n)
\subset
\mbox{gr}^{i+k} {\cal{F}}_{\cal{Z}} (\pi_1, \dots, \pi_n).$$ 

Thus $\gr {\cal{F}}_{\cal{Z}} (\pi_1, \dots, \pi_n)$ is also
a finite--dimensional graded cyclic $\toki{a}$--module.

\begin{statement}{Definition}
Let $\pi_1, \dots, \pi_n$ be finite--dimensional
graded $\toki{a}$--modules with cyclic vectors $v_1, \dots, v_n$,
${\cal{Z}} \in \conf$. 
Then the {\rm fusion product}
$(\pi_1 * \dots *\pi_n)({\cal{Z}})$ is the graded  $\toki{a}$--module
$ \gr {\cal{F}}_{\cal{Z}} (\pi_1, \dots, \pi_n)$ with cyclic vector $v_1
\otimes\dots\otimes v_n \in \mbox{gr}^0 {\cal{F}}_{\cal{Z}}$.
\end{statement}

\begin{statement}{Conjecture}\label{gen}
Let $\gtc{g}$ be a simple Lie algebra, $\pi_1,
\dots, \pi_n$ be cyclic graded $\toki{g}$--modules. Then

\rf For any ${\cal{Z}}_1, {\cal{Z}}_2 \in \conf$ we have
$$(\pi_1 * \dots * \pi_n) ({\cal{Z}}_1) \cong (\pi_1 * \dots * \pi_n)
({\cal{Z}}_2)$$
as $\toki{g}$--modules. 

\rr Fusion product is associative up to isomorphism. It means that
$$ ((\pi_1 * \pi_2) * \pi_3) \cong (\pi_1 * (\pi_2 * \pi_3)) \cong (\pi_1
* \pi_2 * \pi_3)$$
and similar isomorphisms for $4$--fold products and further.
\end{statement}

We prove some cases of this conjecture in the section \ref{count}. 

Of course, (i)
implies that the filtered $\gtc{g}$--module ${\cal{F}}_{\cal{Z}} (\pi_1,
\dots, \pi_n)$ doesn't depend on ${\cal{Z}}$. Nevertheless, the structure
of
$\toki{g}$ module on it depends on ${\cal{Z}}$. Now let us discuss the
relation between these two constructions.

Fusion product can be considered as the extension of the filtered tensor
product from $\conf$ to the origin. Let $\pi_1, \dots, \pi_n$ be
$\gtc{g}$--modules. 
For ${\cal{Z}} \in \Complex^n$
consider the ideal
$$ I({\cal{Z}})  = \gtc{g} \otimes \prod_{i=1}^n (t - z_i) \Complex[t]
\subset \toki{g}.$$

Then we easily obtain

\begin{statement}{Proposition} 

\rf The ideal $I({\cal{Z}})$ annihilates the
filtered tensor product
${\cal{F}}_{\cal{Z}}(\pi_1, \dots, \pi_n)$.

\rr The ideal $I(0) = \gtc{g} \otimes t^n\Complex[t]$  annihilates the
fusion product $(\pi_1 * \dots *
\pi_n)({\cal{Z}})$.
\end{statement}

We can also  extend the filtration to any point ${\cal{Z}} \in \Complex^n$.
Suppose that representations $\pi^1_i, \dots \pi^{s_i}_i$ correspond the
point $z_i \in \Complex$, $i=1\dots m$. Then we should consider distinct
"virtual" coordinates $z^j_i$, $j=1\dots s_i$, and graded
$\toki{g}$--modules
$$\Pi_i = (\pi^1_i * \dots * \pi^{s_i}_i)(z^1_i, \dots z^{s_i}_i).$$
Finally  we can define
$${\cal{F}}_{\cal{Z}}(\pi^j_i) = {\cal{F}}_{(z_1, \dots z_m)} (\Pi_1,
\dots, \Pi_m).$$

Clearly this $\toki{g}$--module is annihilated by the ideal
$I({\cal{Z}})$.

Indeed this is an extension of the filtration to the space larger then
$\Complex^n$. This space can be projected to $\Complex^n$ by forgetting
all "virtual" coordinates. And conjecture~\gref{init}{gen}~(i) implies
that the result depends only on a point of $\Complex^n$. 

Generally, we can consider any superposition of fusion products and
filtered tensor product as the filtered tensor product calculated at a
point ${\cal{Z}}$ of a certain space ${\cal{M}}_n$. This space resembles
the
compactification of the configuration space of $n$ distinct points
constructed in \cite{FMp}.

${\cal{M}}_n$ is stratified. Its strata are enumerated by trees with $n$
input vertices. These trees describe the method of superposition. For
example $\conf \subset {\cal{M}}_n$ as the stratum with $n$ input
vertices
connected to the output vertex.

Any point of ${\cal{M}}_n$ is determined by a tree and complex numbers
attached to each vertex except the output vertex. These complex numbers
are parameters of fusion and filtered tensor products.

Consider the natural map ${\cal{M}}_n \to \Complex^n$ that discards all
the coordinates except those used in the external filtered tensor product.
Then conjecture~\gref{init}{gen} states that for any ${\cal{Z}} \in
{\cal{M}}_n$ the
extended product ${\cal{F}}_{\cal{Z}}$ depends only on the image of
${\cal{Z}}$ in $\Complex^n$.

\section{$\itname{sl}_2$--case and generalizations.}\label{count}\sbros

Here we discuss how to calculate the filtration and the fusion product.
Also we will attempt to clarify the situation discussed in \cite{Sz}. 

\subsection{Computation of the fusion product}\label{calc}

\def \slc {{sl_2^\Complex}}

Let $e,f,h$ be the standard basis in $\sltwo$, and let $e_i, f_i, h_i$ be
a basis in $\sltwo^\Complex$, $e_i = e \otimes t^i,\, f_i = f \otimes
t^i,\, h_i = h \otimes t^i,\, i \ge 0$. Let $\pi_1, \dots, \pi_n$ be
finite--dimensional graded
$\sltwo^\Complex$--modules with cyclic vectors  $v_1, \dots, v_n$.

Suppose that $h_i v_k = 0$, $e_i v_k = 0$ for $i \ge 0$,
$k = 1, \dots, n$. Then each $\pi_k$ can be considered as an
${\frak f}^\Complex$--cyclic module, where $\frak f$ is the
one--dimensional abelian
algebra generated by $\{ f\}$. It is clear that the fusion $(\pi_1 * \pi_2
* \dots *
\pi_n)(z_1, \dots, z_n)$ of ${\frak f}^\Complex$--modules is isomorphic to
the restriction of the fusion $(\pi_1 * \pi_2 * \dots * \pi_n)(z_1, \dots,
z_n)$ of
$\sltwo^\Complex$--modules to ${\frak f}^\Complex$. So it is enough to
calculate the fusion product of ${\frak f}$--modules. But sometimes the
structure of $sl_2$--module is very useful.

Consider in this section the following subcategory $L$ of
finite--dimensional graded
$\sltwo^\Complex$--modules with cyclic vector. Module $(\pi, v)$ is an
object of $L$ if $v$ is killed by
$\{h_i, \, i>0 \}$ and $\{e_i, \, i\ge 0\}$. Morphisms are
$\sltwo^\Complex$--equivariant maps
transforming cyclic vectors to cyclic vectors.
We present the version of Demazure
reflection functor (see \cite{D}) acting on $L$. 

Let $\frak u$ be the subalgebra in
$\sltwo^\Complex$ of codimension one with basis $\{e_i, i \ge 0;\, h_i,
i \ge 0;\, f_i, i \ge 1\}$. Let $S$ be a finite--dimensional $\frak
u$--module. For simplicity suppose that $S$ is $h_0$--semisimple. 
By $\name{Ind}({\frak u}; S; \sltwo^\Complex)$  denote the
$\sltwo^\Complex$--module induced from the $\frak u$--module $S$. 
Then there exists the maximal
finite--dimensional $\sltwo^\Complex$--quotient of $\name{Ind}({\frak u};
S;
\sltwo^\Complex)$. Namely, let $\name{Ind}_{fin}({\frak u}; S;
\sltwo^\Complex)$ be the quotient $\name{Ind}({\frak u}; S;
\sltwo^\Complex) / J$ where $J$
is the intersection of all the submodules such the quotient is
finite--dimensional.

\def \infin {\name{Ind}_{fin}}

\begin{statement}{Proposition}
$\infin(\gtc{u}; S; \slc)$ is a finite--dimensional
$\slc$--module.
\end{statement}

\begin{proof}
Let $\gtc{u}_0$ be the two--dimensional subalgebra of $\slc$ with basis
$\{ ht^0, et^0 \} $. Then, clearly, the vector space $\infin(\gtc{u}; S;
\slc)$ is a
quotient of the finite--dimensional vector space $\infin(\gtc{u}_0; S;
sl_2)$.  
\end{proof}

Now let $S \in L$. Then we have

\begin{statement}{Proposition}\label{id}
$\infin(\gtc{u}; S;
\slc) \cong \infin(\gtc{u}_0; S; sl_2)$ as bi--graded 
$sl_2$--modules. 
\end{statement}

\begin{proof}
Let $\infin(\gtc{u}_0; S; sl_2) = \name{Ind}(\gtc{u}_0; S; sl_2)/J$.
By $v$ denote the cyclic vector of $S$, by $\lambda$ denote the
$h$--weight of $v$.
Then, clearly, $J$ is formed by elements $e^k f^l P(ft, ft^2, \dots) v$,
such that $l > \lambda - 2 \deg (P)$.

As $et^i v = ht^i v = 0$ for $i>0$, it is easy to show that $J$ is an
$\slc$ -- submodule, so we have the statement of the proposition.
\end{proof}

Thus we can investigate properties of this functor using standard facts
about the category ${\cal{O}}$ for $sl_2$.

\begin{statement}{Proposition}\label{idn}
Let $S = \bigoplus S_i$ be a graded $\gtc{u}_0$--module, $h$ acts on
$S_i$ by scalar $i$. By $\pi(i)$
denote the irreducible $sl_2$--module of dimension $i+1$. 

\rf $\infin(\gtc{u}_0; S; sl_2) = \bigoplus_i \pi(i) \otimes (S_i/e^{i+1}
S_{-i-2})$, where $sl_2$ acts only on $\pi(i)$.

\rr $\infin(\gtc{u}_0; S; sl_2)$ is right exact and  $L^i\infin(\gtc{u}_0;
S; sl_2) = 0$ for $i>1$.

\rr $L^1\infin(\gtc{u}_0; S; sl_2) = \bigoplus_i \pi(i)
\otimes \name{Ker}_{S_{-i-2}} (e^{i+1})$, where $sl_2$ acts only on
$\pi(i)$.
\end{statement}

In particular, the Euler characteristic $\ch L^0 - \ch L^1$ of this
functor depends only on the character of the representation $S$. Namely,
by the proposition above it is
equal to  $\sum \, (\dim S_i - \dim S_{-i-2})\cdot \ch \pi(i)$. It
resembles
the situation
in \cite{D}. Now let us continue our considerations.

Let $\theta$ be the homomorphism ${\frak u} \to \sltwo^\Complex$ such that
$\theta(e_i) = e_{i+1}$, $\theta(h_i) = h_i$, $\theta(f_i) = f_{i-1}$.
With each complex number $\lambda$ we associate a $1$--dimensional
representation
$\chi_\lambda: {\frak u} \to \Complex$, such that $\chi_\lambda(h_0) =
\lambda$, $\chi_\lambda(f_i) = 0$, $\chi_\lambda(h_i) = 0$, $i \ge 
1$; $\chi_\lambda(e_j) = 0$, $j \ge 0$.

\begin{statement}{Definition}
Let $\lambda \in \Complex$. Functor ${\cal T}_\lambda$ acts from $L$
to $L$ in the following way. Let $\pi \in L$, then $\theta^*(\pi)$ is 
a representation of $\frak u$. Let
%*
\begin{displaymath}
{\cal T}_\lambda(\pi) = \name{Ind}_{fin}({\frak u}; \theta^*(\pi) \otimes
\chi_\lambda; \sltwo^\Complex).
\end{displaymath}
%*
\end{statement}

Let $\pi_1, \dots, \pi_n$ be finite--dimensional representations of
$\sltwo$ with highest vectors $v_i$ ($e v_i = 0$). Order them such that  
$\dim \pi_1 \le \dim \pi_2 \le \dots \le \dim \pi_n$. Let $\lambda_1 \le
\lambda_2 \le \dots \le \lambda_n$ be their highest weights.

\begin{statement}{Theorem}\label{main}
The $\sltwo^\Complex$--representation $(\pi_1 * \pi_2 * \dots *
\pi_n)(z_1,
\dots, z_n)$ (recall that the points $z_1, \dots, z_n$ are pairwise
distinct) is isomorphic to the representation
%*
\begin{displaymath}
({\cal T}_{\lambda_n} \dots ({\cal T}_{\lambda_2}({\cal
T}_{\lambda_1}(1)))
\dots),
\end{displaymath}
%*
where $1$ is trivial representation. In particular, the
$\sltwo^\Complex$--module $(\pi_1 * \pi_2 * \dots * \pi_n)(z_1, \dots,
z_n)$ 
does not depend on the choice of the configuration $z_1, \dots, z_n$.
\end{statement}

\begin{proof} Let $v = v_1 \otimes v_2 \otimes \dots \otimes v_n$ be the  
cyclic vector in
$\pi_1(z_1) \otimes \dots \otimes \pi_n(z_n)$. 

First of all, remark that as $\pi_i$ 
are cyclic $\gtc{f}$--modules, their fusion product and filtered tensor
product are also cyclic $\toki{f}$--modules or, in other words, quotients
of
$\Complex [f_0, f_1, \dots] v$
with the natural action of $sl_2^\Complex$.

Let ${\cal{Z}} = (z_1, \dots, z_n)$, $z\in \Complex$. By ${\cal{Z}} + z$
denote the shifted vector $(z_1 + z, \dots, z_n + z)$. Clearly,
$${\cal{F}}_{{\cal{Z}} + z} (\pi_1,\dots, \pi_n) \cong
{\cal{F}}_{\cal{Z}}
(\pi_1,\dots, \pi_n) (z),$$
hence $(\pi_1 * \dots * \pi_n)({\cal{Z}} + z) \cong (\pi_1 * \dots *
\pi_n)({\cal{Z}})$ and we can suppose that in our case $z_n=0$.

Then the elements $f_i$ for
$i\ge 1$ act trivially on $v_n$, so
$$\Complex[f_1, f_2,\dots](v_1 \otimes \dots \otimes v_n) \subset
\pi_1\otimes \pi_2 \otimes \dots \otimes \pi_{n-1} \otimes v_n.$$   
We have now the picture where only
the representations $\pi_1(z_1), \dots, \pi_{m-1}(z_{m-1})$ are involved.

Consider the fusion product $(\pi_1 * \dots *\pi_{n-1})({\cal{Z}})$,
denote by $u$ its highest vector. We know
that this representation is a quotient of $\Complex [f_0, f_1\dots] u$
by a certain
ideal $J$. We can apply the map $\theta : \Complex[f_1, f_2,\dots] \to
\Complex [f_0, f_1, \dots]$ to this ideal.

\begin{Statement}{Lemma}
$\theta (J) \subset J$.
\end{Statement}

\begin{proof} 
As $\pi_i \cong \Complex[f]v_i$ one can consider the
automorphism of $\pi_i$ mapping $f^k v_i$ to $z_i^k f^k v_i$. Applying
the product of these automorphisms to the $sl_2^\Complex$--module
$\pi_1(z_1) \otimes
\dots \otimes \pi_n(z_n)$ we substitute the action of $f_i$ by the action
of $f_{i+1}$. But the module $(\pi_1 * \dots * \pi_{n-1})({\cal{Z}}) =
({\cal T}_{\lambda_{n-1}} \dots ({\cal T}_{\lambda_2}({\cal
T}_{\lambda_1}(1)))\dots)$ 
does not change under this automorphism. Thus
$f_1, f_2, \dots$ must satisfy the same relations as $f_0, f_1, \dots$.
\end{proof}  

Suppose now that the statement of the theorem holds for
$n-1$ representations
$\pi_1(z_1), \dots, \pi_{n-1}(z_{n-1})$. Let
$$N = \Complex[f_1, f_2, \dots] / \theta(J) = \theta^* ((\pi_1  * \dots * 
\pi_{n-1})(z_1, \dots, z_{n-1})).$$ 
Due to the lemma above we have 
the surjection $N \stackrel{\mbox{\scriptsize \(\phi\)}}{\to}
\Complex[f_1, f_2, \dots] v$, where $v$ is the vacuum vector in $(\pi_1 *
\dots * \pi_n)(z_1, \dots, z_n)$.  Note that both $N$ and $\Complex[f_1,
f_2, \dots] v$ have natural structures of ${\frak u}$--modules, and $\phi$
is compatible with these structures. Therefore we have the map
$$\name{Ind}_{fin} ({\frak u}; N; \sltwo^\Complex) \to \name{Ind}_{fin}
({\frak u}; \Complex[f_1, f_2, \dots] v; \sltwo^\Complex).$$ 
By proposition~\gref{count}{idn} (ii) it is also a
surjection. Then, as the fusion product is finite--dimensional, we have
the surjection
$$\name{Ind}_{fin} ({\frak u}; N; \sltwo^\Complex) \to \Complex[f_0,
f_1, \dots] v = (\pi_1 * \dots * \pi_n)(z_1, \dots, z_n).$$

Thus it remains to 
prove that $\dim(\name{Ind}_{fin} ({\frak u}; N; \sltwo^\Complex))$ is
equal to the dimension of $\pi_1 \otimes \pi_2 \otimes \dots \otimes
\pi_n$.

\begin{Statement}{Lemma} $L^1 {\cal{T}}_{\lambda_n}( \pi_1* \dots 
* \pi_{n-1}) = 0$
\end{Statement}

\begin{proof}
Due to the proposition \gref{count}{idn} we need to prove that for any $k$
$(et^1)^k$ is injective on the subspace of $\pi_1 * \dots * \pi_{n-1}$
of $h$--weight $-(k+1)/2$. It is equivalent to the fact that 
$(ft^1)^k$ is injective
on the
subspace of $h$--weight $(k+1)/2$. But this fact can be easily checked by
induction using the condition that $\lambda_1 \le \lambda_2 \le \dots \le
\lambda_n$.
\end{proof}

Due to the propositions \gref{count}{id} and \gref{count}{idn} we know
the dimension of ${\cal{T}}_{\lambda_n}( \pi_1 * \dots 
* \pi_{n-1})$ and can compare it with $\dim (\pi_1 \otimes \dots \otimes
\pi_n)$.
\end{proof}

So the $sl_2^\Complex$--module $(\pi_1 * \dots * \pi_n)(z_1, \dots, z_n)$ 
does not depend on $z_1, \dots, z_n$. Therefore we
can write $\pi_1 * \dots * \pi_n$ without $\{z_i\}$. We see that the
$\Complex[f_0, f_1, \dots]$--module $\pi_1 * \dots * \pi_n$ is a quotient
$\Complex[f_0, f_1, \dots]/J$. The algebra $\Complex[f_0, f_1, \dots]$ has
a bi--grading when the degree of the generator $f_i$ is equal to $(1,i)$.
$J$ is a bi--graded ideal. Then
%*
\begin{displaymath}
\pi_1 * \dots * \pi_n = \bigoplus_{r,s} V_{r,s}
\end{displaymath}
%*
where $V_{r,s}$ is the part of degree $(r,s)$.

Introduce the character
%*
\begin{displaymath}
\name{ch}(\pi_1, \dots, \pi_n; z,q) = \sum (\dim V_{r,s}) \cdot z^r \cdot
q^s.
\end{displaymath}
%*

As a consequence of the theorem above we can write down the following
recurrent relation. Let $\pi_1, \dots, \pi_n$ be irreducible
$sl_2$--modules such that $\dim \pi_n \ge \dim \pi_j$ for all $j$. Then
%*
\begin{displaymath}
\name{ch}(\pi_1, \dots, \pi_n; z,q) = \name{ch}(\pi_1, \dots, \pi_{n-1};
qz,q) + \name{ch}(\pi_1, \dots, \pi_{n-1}, \widetilde{\pi_n}; z,q) \cdot z
\end{displaymath}
%*
where $\widetilde{\pi_n}$ is a $\sltwo$--irrep such that 
$\dim \widetilde{\pi_n} = \dim \pi_n - 1$.

\begin{Statement}{Remark}
We can similarly prove some associativity relations, for example
that if $\dim \pi_3 \ge \dim \pi_2, \dim \pi_1$ then
$$ ((\pi_1 * \pi_2) * \pi_3) \cong (\pi_1 * \pi_2 * \pi_3).$$
Proof of other associativity relations is more subtle and
we are going to present it later. 
\end{Statement}

\subsection{Connection with the representation theory of the current
algebra}\label{con}

\def \slhat {\widehat{\sltwo}}
\def \Integer {{\Bbb Z}}

Let $\slhat$ be the affine Kac--Moody algebra associated with $\sltwo$. Fix
a non-negative integer $k$. It is known that there are $k+1$ irreducible
integrable representations of $\slhat$ at level $k$. Denote them by
$L_0, L_1, \dots, L_k$. Here $L_j$ is the representation with the highest
weight $(j,k-j)$. Consider the following basis in $\tsl$. Let 
$e_i = e \otimes t^i$, $h_i = h \otimes t^i$, $f_i = f
\otimes t^i$, $K$ be the central element and $d$ be the energy element.
Then the vacuum vector in $L_j$ is killed by $e_i$,
$i \ge 0$, $h_i$, $i > 0$, and $f_i$, $i > 0$.

Here we need a standard fact from representation theory of
current algebras.

\begin{statement}{Proposition}\label{rel} Consider the operators
$\{f_i\}$ on $L_j$. They satisfy only the relation $f^{k+1}(z) =
0$, where $f(z) = \sum f_i \cdot z^{-i-1}$. It means that
%*
\begin{displaymath}
\sum_{{i_p \in \bbb{Z}}\atop{i_1 + i_2 + \dots + i_{k+1} = s}}
f_{i_1} f_{i_2} \dots f_{i_{k+1}} = 0
\end{displaymath}
for any $s$.

\end{statement}

In particular, if $v$ is the vacuum vector then

\begin{displaymath}
\sum_{{i_p \in \bbb{Z}, i_p \le 0}\atop{i_1 + i_2 + \dots + i_{k+1} = s}}
f_{i_1} f_{i_2} \dots f_{i_{k+1}} v = 0
\end{displaymath}
for any $s$.

Let $v_j$ be the vacuum vector in $L_j$ and let $W_j \subset L_j$ be
the subspace $\Complex
[f_0, f_{-1}, f_{-2}, \dots] \cdot v_j$. Then one can show that
$W_j$ is isomorphic to the quotient
$\Complex [f_0, f_{-1}, f_{-2}, \dots] / J$ where $J$ is the ideal
generated by the elements
%*
\begin{displaymath}
f_0^{j+1} = 0\quad \mbox{and} \
\sum_{{i_p \in \bbb{Z}, i_p \le 0}\atop{i_1 + i_2 + \dots + i_{k+1} = s}}
f_{i_1} f_{i_2} \dots f_{i_{k+1}} = 0, \quad s=0,-1,-2,\dots.
\end{displaymath}
%*

Let $W_j(N) = \Complex [f_0,
f_{-1}, \dots, f_{-N}] \cdot v_j$ be a subspace of $L_j$, $N \ge 0$. Note
that the space $W_j(N)$
is stable with respect to the action of operators $e_i,\, i \ge N$,
$h_j,\, j \ge 0$, and $f_k,\, k \ge -N$. Clearly, the Lie algebra
formed by $\{e_i, h_j, f_k;\, i \ge N, j \ge 0, k \ge -N\}$ is
isomorphic to $\sltwo^\Complex$.

\begin{statement}{Proposition}\label{str}
The $\sltwo^\Complex$--module $W_j(N)$ is isomorphic to the fusion
product $[j] * [k] * \dots * [k]$, where $[r] = S^r [1]$ is
the irreducible representation of $\sltwo$ of dimension $r+1$, $[k]$
appears $N$ times.
\end{statement}

\begin{proof}
Induction on $N$. The fact that ${\cal{T}}_k({\cal{T}}_k( \dots
\mbox{($N-1$
times)}
\dots  {\cal{T}}_j(1)))
= W_j(N-1)$ implies that there exists a surjection 
$${\cal{T}}_k({\cal{T}}_k( \dots \mbox{($N$ times)} \dots {\cal{T}}_j(1)))
\to W_j(N).$$
The rest of proof is comparing the dimensions.
\end{proof}

Consider now the vacuum representation $L_0$. Let $\{w(j)\}$ be the set of 
extremal vectors in $L_0$, $j \in \Integer$. Each vector $w(j)$ is defined
up to a scalar by the properties that $e_i w(j) = 0$, $i \ge -2j$;
$h_i w(j) = 0$, $i > 0$;
$f_i w(j) = 0$, $i \ge 2j$.
Consider the standard inclusion $\sltwo^\Complex \hookrightarrow \slhat$,
where the
image is $\{e_i, h_i, f_i,\, i \ge 0\}$ (it corresponds
to the case $N = 0$ in the notation of proposition~\gref{count}{rel}). Let
$R(j)
= U(\sltwo^\Complex) w(j),\, j \ge 0$.

We know that $R(j) = ([k] * \dots * [k])$. Then
conjecture~\gref{init}{gen}
implies that 
$\bigl(R(j_1) * R(j_2)\bigr)(z_1, z_2)$ is isomorphic to $R(j_1 + j_2)$
for any $z_1 \ne z_2$.

\begin{Statement}{Remark}
The vacuum representation $L_0$ is the direct limit of the following
sequence
of inclusions: $\Complex = R(0) \hookrightarrow R(1) \hookrightarrow R(2)
\hookrightarrow \dots$. Therefore, $L_0$ is the direct limit $\Complex
\hookrightarrow [k] * [k] \hookrightarrow [k] * [k]
* [k] * [k] \mbox{\space (\(4\) times)} \hookrightarrow \dots$.
Generally, the representation $L_j$ can be constructed as the direct
limit of the sequence of $\sltwo^\Complex$--homomorphisms
%*
\begin{displaymath}
[j] \hookrightarrow [j] * [k] * [k] \hookrightarrow
[j] * [k] * [k] * [k] * [k] \hookrightarrow
\dots
\end{displaymath}
%*
The inclusion $[j] * [k]^{2s-2} \hookrightarrow [j] *
[k]^{2s}$ is defined in the following way. Cyclic vector from
$[j] * [k]^{2s-2}$ goes to $(f_{2s})^j (f_{2s-1})^{k-j} \cdot
v$, where $v$ is the cyclic vector in $[j] * [k]^{2s}$
\end{Statement}

\def \of {{\overline{f}}}

We can easily apply this consideration to computation of the fusion
product. 
Let $\pi_1, \dots, \pi_N$ be $sl_2$--modules. Then the
$\gtc{f}^\Complex$--module $\pi_1 * \dots
* \pi_N$ is a quotient of $\Complex [f_0, \dots, f_{N-1}]$.
It seems wise to rename the variables. Namely, let $\of_0 = f_{N-1}$,
$\of_1 = f_{N-2}$, \dots, $\of_{N-1} = f_0$.

Suppose all the representations are $[k]$. Then propositions
\gref{count}{rel} and \gref{count}{str} imply that $[k] * \dots * [k]$
is the image of $\Complex[\of_0, \dots, \of_{N-1}]$ in the quotient of
$\Complex[\of_0, \of_1, \dots]$ by the ideal generated by
$$ \psi_{k+1}(s) =
\sum_{{i_p \in \bbb{Z}, i_p \ge 0}\atop{i_1 + i_2 + \dots + i_{k+1} = s}}
\of_{i_1} \of_{i_2} \dots \of_{i_{k+1}} = 0$$
for $s = 0,1,\dots$.

For $[j] * [k] * \dots * [k]$ we have the similar construction with
the ideal generated by $\psi_{k+1}(s)$ for all $s$ and $\psi_{j+1}(0)$.

Probably, there exists a similar construction for arbitrary
representations $\pi_1, \dots, \pi_n$. Let $k+1$ be the maximal dimension
of the representations.
Then, probably, in this case the ideal is generated by elements
$\psi_{k+1}(s)$ for all $s$
and elements $\psi_{r}(s)$ for $r<k+1$ and small $s$.

\subsection{Relation to Schubert varieties}\label{var}

\def \Proj {{\Bbb P}}

Let $\gtc{a}$ be a Lie algebra,
$\pi$ be a cyclic finite--dimensional 
${\frak a}^\Complex$--module with cyclic vector $v$. Let ${\frak
a}^\Complex[\pi]$ be the image of ${\frak a}^\Complex$ in the Lie
algebra $\name{End}(\pi)$, and let ${\frak A}[\pi]$ be the group
corresponding to ${\frak a}^\Complex[\pi]$.
Consider the projectivization $\Proj(\pi)$ of the vector space $\pi$. The
group $\gtc{A}[\pi]$ acts on $\Proj(\pi)$. 
By $O(\pi)$ we denote the closure of the $\gtc{A}[\pi]$--orbit
of the point $\Complex \cdot v$ in $\Proj(\pi)$. 
We have the canonical bundle $\xi$ on $O(\pi)$ such that
in reasonable cases $H^0(\xi) = \pi$.

\begin{statement}{Proposition}
Suppose $\frak a$ be $\sltwo$, $\pi$ be a
finite--dimensional irreducible representation. Then the manifold $O(\pi *
\dots * \pi)$ is one of the Schubert varieties on the affine flag variety
for the Lie algebra $\slhat$ (see \cite{PS}). By Schubert variety here we
mean the closure of an orbit of
the group $\name{SL}_2(\Complex[t])$ on the affine flag variety.
\end{statement}

Actually this statement is a reformulation of the results of
section~\ref{con}. In section~\ref{con} we
described a construction of the integrable representations as the
direct limit of the ``fusion'' tensor products. Standard geometric
approach realize the irreducible representation in the space of sections of
the line bundle on the affine flag manifold. Consider for example the
vacuum representation $L_0$ (at level $k$). $L_0$ can
be realized in the space of sections of the line bundle $\zeta^k$ on the
``grassmannian''
$\name{SL}_2(\Complex[t,t^{-1}])/\name{SL}_2(\Complex[t])$.
Schubert varieties on the grassmannian are labelled by non--negative
integers. Let us
denote
them $\name{Sh}(2j)$, $\dim \name{Sh}(2j) = 2j$. The Lie algebra
$\sltwo^\Complex$ acts on $\name{Sh}(2j)$, therefore $\sltwo^\Complex$
acts on the space $H^0(\name{Sh}(2j), \xi^k)$. Here we denote by $\xi$
the restriction of $\zeta$ to $\name{Sh}(2j)$.

\begin{statement}{Proposition}
$O(\pi * \dots * \pi) = \name{Sh}(2j)$ ($\pi$ appears $2j$ times).
The $\sltwo^\Complex$--module $H^0(\name{Sh}(2j), \xi^k)$ is dual to the
representation $([k] * \dots * [k])$ ($2j$ times).
\end{statement}

To investigate the product of odd number of representations consider the
``twisted'' $\slhat$--grassmannian. It is the
quotient 
$\name{SL}_2(\Complex[t,t^{-1}])/G$
where $G$ is the group with Lie algebra 
$$\left(\begin{array}{cc}a(t) & tb(t) \\ t^{-1} c(t) &
-a(t)\end{array}\right),$$
where $a,b,c$ are polynomials in $t$. The orbits
of $\name{SL}(\Complex[t])$ now have odd dimension, and we denote them
by $\name{Sh}(2j+1)$. The  
$\sltwo^\Complex$--module $H^0(\name{Sh}(2j+1), \xi^k)$ is dual to the
fusion product $[k] * \dots * [k]$
($2j+1$ times). So $O(\pi * \dots * \pi) = \name{Sh}(2j+1)$, $\pi$
appears $2j+1$--times.

Consider now the case of different representations. Suppose, we have 
 $N$ irreducible representations $\pi_1, \pi_2, \dots,
\pi_N$ of $sl_2$, where $r_1$ are of dimension $k_1$, $r_2$ are of
dimension $k_2$,
\dots, $r_s$ are of dimension $k_s$, $k_1 < k_2 < \dots < k_s$, and $r_1 +
r_2
+ \dots + r_s = N$.

\begin{statement}{Proposition}
There is a surjection $O(\pi_1 * \pi_2 * \dots * \pi_N) \to \name{Sh}(N)$.
\end{statement}

\begin{proof}
In section~\ref{calc} we show that there is a map of $\Complex[e_0,
e_1, \dots]$--modules
$\pi_1 * \pi_2 * \dots * \pi_N \to \pi_1 * \pi_2 * \dots *
\widetilde{\pi_N}$. Here we suppose that $\dim \pi_N = k_s$ (the biggest)
and $\dim \widetilde{\pi_N} = k_s - 1$. Suppose that $\dim \pi_1 = k_1$
(the smallest). Iterating this procedure, we get a $\Complex[e_0, e_1,
\dots]$--homomorphism $\pi_1 * \pi_2 * \dots * \pi_N \to \pi_1 * \dots *
\pi_1$ ($N$ times). Note that 
the cyclic vector $v$
in $\pi_1 * \pi_2 * \dots * \pi_N$ 
is killed by $\{f_i, h_i\}$. We deduce from this, that $O(\pi_1 *
\pi_2 * \dots * \pi_N)$ is the closure of the orbit of $\Complex\cdot v$
with
respect to the action of the group with the Lie algebra $\{e_i\}$.
Therefore we have a well--defined map of the orbit on $\Proj(\pi_1 *
\pi_2 * \dots * \pi_N)$ to the orbit on $\Proj(\pi_1 * \dots * \pi_1)$.
The proof that
this map can be extended to the closure is more subtle. We will present it
in another place.
\end{proof}

\begin{statement}{Proposition}
There is a fibration $O(\pi_1 * \pi_2 * \dots * \pi_N) \to
\name{Sh}(r_s)$.
The fiber is isomorphic to $O(\pi_1 * \dots * \pi_{N-r_s})$.
\end{statement}

Note that $\name{Sh}(1) = \Complex P^1$. As a consequence we have

\begin{statement}{Proposition}
Suppose that $\dim \pi_1 < \dim \pi_2 < \dots < \dim \pi_N$. Then $O(\pi_1
* \pi_2 * \dots * \pi_N)$ is a smooth manifold. Therefore the map $O(\pi_1
* \pi_2 * \dots * \pi_N) \to \name{Sh}(N)$ gives us a resolution of
singularities of the manifold $\name{Sh}(N)$.
\end{statement}

\subsection{Case of abelian Lie algebra; $\itname{sl}_{n+1}$--case}

\def \A {{\frak a}}

Let $\A$ be an abelian algebra. We saw that the $\sltwo$--case is actually
deeply related to the case of abelian $1$--dimensional algebra. More
precisely, let $x \in \A$, $\dim \A = 1$, and let $M(r)$ be the
representation
of $\A$ in the space $\Complex[x]/x^r \Complex[x]$. The space $M(r)$ can be
treated as a representation of $\A^\Complex$. The study of the
$\A^\Complex$--module $(M(r_1) * M(r_2) * \dots * M(r_k))(z_1, \dots,
z_k)$
is then a part of the corresponding problem for $\pi_1 * \dots * \pi_k$
where $\pi_j$ are finite--dimensional representations of $\sltwo$.

We consider now a generalization of this $1$--dimensional case. Let $\A$
be
an $n$--dimensional algebra, and $x_1, \dots, x_n$ be a basis in $\A$.
Denote by $S^r$ the quotient 
$\Complex[x_1, \dots, x_n]/J_r$,
where $J_r$ is spanned by monomials $x_1^{i_1} x_2^{i_2} \dots x_n^{i_n}$
such that $i_1 + i_2 + \dots + i_n > r$. Surely, $S^r$ is a cyclic
$\A$--module and
$1$ is its cyclic vector.

Consider now the Lie algebra $\itname{sl}_{n+1}$, and let $V$ be the
canonical irreducible $(n+1)$--dimensional representation of
$\itname{sl}_{n+1}$.
Fix also the standard decomposition $\itname{sl}_{n+1} = W \oplus
\itname{gl}_n
\oplus W^*$, $\dim W = n$, $W$ and $W^*$ are abelian subalgebras,
$\itname{gl}_n$ acts on $W$ by the standard way. Let $S^r(V)$ be the
$r$--th symmetric power of $V$. $S^r(V)$ is an irreducible representation
of
$\itname{sl}_{n+1}$. Consider the vector $v$ in $S^r(V)$ such
that
$W^* \cdot v = 0$ and $\Complex \cdot v$ is the one--dimensional
representation
of $\itname{gl}_n \subset \itname{sl}_{n+1}$. Clearly, it is a cyclic
vector.

\begin{statement}{Proposition}
$S^r(V)$ is isomorphic to $S^r$ as cyclic $W$--modules.
\end{statement}

\begin{statement}{Theorem}

\rf The $\itname{sl}_{n+1}$--module $(S^{r_1}(V) * S^{r_2}(V) * \dots *
S^{r_k}(V))(z_1, \dots, z_k)$ does not depend on the choice of the set
$\{z_1, \dots, z_k\}$.

\rr The restriction of $S^{r_1}(V) * S^{r_2}(V) * \dots * S^{r_k}(V)$
to the abelian algebra $W^\Complex$ is isomorphic to $S^{r_1} * S^{r_2} *
\dots * S^{r_k}$. In particular, for abelian algebra $W$ the
$W^\Complex$--module $(S^{r_1} * S^{r_2} * \dots * S^{r_k})(z_1, \dots,
z_k)$ does not depend on the choice of configuration $\{z_1, \dots,
z_k\}$.
\end{statement}

The proof is basically the same as in the $\sltwo$--case. It uses a
version
of the Demazure reflections.

\subsection{ Geometric picture for the $\itname{sl}_{n+1}$.}

For $\itname{sl}_{n+1}$ there is a collection of manifolds
$\name{sh}(\lambda)$
labelled by dominant weights $\{\lambda\}$ of $\itname{sl}_{n+1}$.
Lie algebra $\itname{sl}_{n+1}$ acts on $\name{sh}(\lambda)$. These
manifolds are Schubert varieties on the ''grassmannians'' (see \cite{PS}).
Each
$\name{sh}(\lambda)$ is equipped with a linear bundle $\xi$. If $\lambda$
is a
fundamental weight then $H^0(\name{sh}(\lambda),\xi)$ is isomorphic to
$\pi(\lambda)$,
the corresponding fundamental representation of $\itname{sl}_{n+1}$.

Fix in each irreducible representation the highest
vector as a cyclic vector. Then for the set of irreducible
representations $\pi_1, \dots, \pi_n$ we can construct the fusion product
$\pi_1 * \dots * \pi_n$ 
and define $O(\pi_1 * \dots * \pi_n)$ as the closure of the orbit of
the highest vector in $\Proj(\pi_1 * \dots * \pi_n)$.

Let $\lambda_1, \lambda_2, \dots, \lambda_n$ be fundamental weights of
$\itname{sl}_{n+1}$.

\begin{statement}{Proposition}
Let $m_1, m_2, \dots, m_n$ be a set of non--negative numbers. Suppose we
have
$m_1$ copies of the representation $\pi(\lambda_1 \cdot k)$, $m_2$ copies
of the
representation $\pi(\lambda_2 \cdot k)$, \dots, $m_n$ copies of the
representation $\pi(\lambda_n \cdot k)$, $k \in \Integer,\, k > 0$. Then
$O(\dots)$ (the fusion of this set of representations) is
$\name{sh}(\lambda),\,
\lambda= \sum m_i \lambda_i$. The $\itname{sl}_{n+1}^\Complex$--module
$H^0(\name{sh}(\lambda), \xi^k)$ is isomorphic to
%*
\begin{displaymath}
(\pi(\lambda_1 \cdot k) * \dots * \pi(\lambda_1 \cdot k)) \mbox{\rm
(\(m_1\) times)} *
\dots * (\pi(\lambda_n \cdot k) * \dots * \pi(\lambda_n \cdot k))
\mbox{\rm (\(m_n\)
times)}.
\end{displaymath}
%*
\end{statement}

\begin{Statement}{Remark}
Our construction can be used for resolution of singularities of the
manifold $\name{sh}(\lambda)$. For example, the manifold $N =
O(\pi(\lambda_1
\cdot
k_1) * \dots * \pi(\lambda_m \cdot k_m))$ is non--singular if $k_1 < k_2 <
\dots
<
k_m$, and there is a map $N \to \name{sh}(\sum \lambda_i)$.
\end{Statement}

\section{Application to Verlinde rule}\label{appl}\sbros

 Here we propose our version of $q$--Verlinde rule. See \cite{FLOT} for
another approach to this problem. 

\subsection{Problem and notation}

Let $\gtc{g}$ be a simple Lie algebra, $\gtc{h}\subset \gtc{g}$ and
$\np \subset \gtc{g}$ be Cartan and nilpotent Borel subalgebras
respectively.  

Let $\tg$ be the corresponding affine algebra. Let us introduce some
notation. 

On the one hand, it is the algebra of $\gtc{g}$--valued algebraic
currents on the circle equipped with the central extension and the energy
element. On
the other hand, it is a Kac-Moody algebra. The standard bases of its
Cartan subalgebra and the dual space are following:
$$ \htg = \left<\{ \chal_0, \chal_1, \dots, \chal_n,
d \}\right> \qquad \htg^* = \left<\{ \alpha_0, \alpha_1, \dots, \alpha_n,
\Lambda_0 \}\right>$$
$$ \left< \chal_i, \alpha_j \right> = a_{ij}, \quad \left< \chal_i,
\Lambda_0 \right> =
\delta_{i,0}, \quad \left< d, \alpha_i \right> = \delta_{i,0}, \quad
\left< d,\Lambda_0 \right> = 0,$$ 
where $\left( a_{ij}\right)$, $i,j = 0\dots n$ is the Cartan matrix of
$\tg$. Note that
the restriction
of this matrix $\left( a_{ij}\right)$, $i,j=1\dots n$ coincides with the
Cartan matrix
of $\gtc{g}$.

Also we need notation for the dual
central element $\delta\in \htg^*$ such that $\left< h, \delta\right> = 0$
for any $h
\in \htg$, 
and notation for the "canonical class" $\rho \in \htg^*$ such that $\left<
\chal_i,
\rho\right>
= 1$ for $i=0\dots n$. These properties define $\delta$ up to
scalar multiplication and $\rho$ up to $\bbb{C}\cdot\delta$.  
For this kind of affine Lie algebras one can take 
$\delta = \alpha_0 + \alpha_M$, where $\alpha_M$ is the maximal root of
$\gtc{g}$.

Now consider the integrable $\tg$--module with the highest weight
$\gamma$.
Then $\gamma$ is a dominant integral weight. Note that the representation
does not
change if we add $c\delta$ to $\gamma$. So we can eliminate the
coefficient before $\alpha_0$ and consider
$$\gamma = \lambda + k\Lambda_0, \quad \lambda = a_1\alpha_1 + \dots
a_n\alpha_n.$$
As $\gamma$ is a dominant integral weight, $\lambda$ is a dominant
integral
weight
for $\gtc{g}$, $k \in \bbb{Z}$ and $\left< \chal_0, \lambda\right> \ge
-k$. Denote
such module by $L_{k, \lambda}$. It is usually called a representation at
level $k$.

Consider the following simplex in $\gtc{h}^*$.
$$\Delta_k = \{ \eta \in \htg^* \mid \ \left< \chal_0, \eta\right> >
-k-1;\
\left< \chal_i, \eta\right> > -1 \}.$$
Then representations at level $k$ are in natural correspondence with
integral weights in $\Delta_k$.

At last, by $W$ denote the Weyl group of $\gtc{g}$ and by $W_{aff}$ denote
the Weyl group of $\tg$. By $w\cdot \gamma$ denote the standard linear
action of Weyl groops, by $w\circ \gamma = w\cdot (\gamma + \rho) - \rho$
denote the "geometric", affine action.

Now we are ready to discuss the problem.

Let $p_1, \dots, p_n$ be distinct complex numbers, ${\cal{P}} =
(p_1, \dots p_n)$. Let $\gtc{a}\subset \gtc{g}$ be a subalgebra.
Consider the
subalgebra $\gtc{g}({\cal{P}}, \gtc{a}) \subset \tg $ of $\gtc{g}$--valued
polynomials $Q$ such that $Q(p_i) \in \gtc{a}$.

To simplify the notation in this section let us "settle" the
representations
$L_{k, \lambda}$ in the infinity. It means that  we substitute
the action of elements $\gtc{g} \otimes t^k$ by the action
of elements $\gtc{g} \otimes t^{-k}$. In particular, elements
$\gtc{g} \otimes t^{-1} \Complex [t^{-1}]$ now act trivially on
the hightest vector.

Introduce the space of coinvariants
$$V(k, \lambda, {\cal{P}}, \gtc{a}) =
L_{k,\lambda}/\gtc{g}({\cal{P}}, \gtc{a}).$$

Indeed in this section either $\gtc{a} = 0$ or $\gtc{a} = \np$.
First case is more usual, in this case the space of coinvariants
has the natural structure of $\gtc{g}$--module. For shortness let
$V(k, \lambda, {\cal{P}}) = V(k, \lambda, {\cal{P}},0)$.
Second case ($\gtc{a} = \np$) is more suitable for calculations.
Here the space of coinvariants has not structure of $\gtc{g}$--module, 
but it admits the action
of $\gtc{h} \subset \gtc{g}$, and, therefore, $\gtc{h}$--grading.

There exists a $\gtc{h}$--equivariant filtration on these spaces.
Namely, $L_{k, \lambda}$ is graded "by $t$" due to the action of $d \in
\htg$, so we have the
induced filtration on the space of coinvariants.

This filtration describes the behavior of coinvariants when ${\cal{P}}$
tends to the origin $0 = (0, \dots, 0)$. We can consider the space $V(k,
\lambda, 0, \gtc{a})$. By general considerations we have
$$\dim V(k, \lambda,0,\gtc{a}) \ge \dim V(k, \lambda, {\cal{P}},\gtc{a})$$
and there exists a subspace $V^c (k, \lambda,0,\gtc{a})
\subseteq V(k,\lambda,0,\gtc{a})$
of
coinvariants that can be deformed to coinvariants at generic
point.

Spaces $V(k, \lambda, 0,\gtc{a})$ and $V^c (k,\lambda,0,\gtc{a})$
are naturally graded "by
$t$" and our filtration is the extension of the grading on
$V^c(k, \lambda,0,\gtc{a})$.

We are convinced that for $sl_n$ we have
$V^c (k, \lambda,0,\np) = V (k,\lambda,0,\np)$ 
and 
$V^c (k, \lambda,0) = V (k,\lambda,0)$, 
but
it can be shown that for the Lie algebra $E_8$ it is not so.

Let
$$ \ch_q (k, \lambda, {\cal{P}},\gtc{a}) = \sum_i q^i \ch (\mbox{gr}^i
 V(k,\lambda,{\cal{P}}),\gtc{a})$$
be the Hilbert polynomial. Our goal is to calculate this polynomial using
the notion of filtered tensor product.

Now let us show that the dimension of the space of coinvariants
doesn't depend on ${\cal{P}}$ and
can be calculated by the Verlinde rule.

\subsection{Verlinde algebra and Verlinde rule}

First of all, recall the definition of the Verlinde algebra at
level $k$.

Let $\ver{}$ be the tensor algebra of virtual
finite--dimensional
$\gtc{g}$--modules. It means that any element of $\ver{}$
is a finite linear combination of basis elements $[\pi_\mu]$
corresponding to irreducible representations $\pi_\mu$.
For any finite--dimensional representation $\pi = \oplus \pi_\nu$ we can
consider the element $[\pi] = \sum [\pi_\nu] \in \ver{}$.
Then the product of two
basis elements $[\pi_\mu]\cdot [\pi_\eta]$ is $[\pi_\mu \otimes
\pi_\eta]$.

Let us define the Verlinde algebra $\ver{k}$ as a quotient of
$\ver{}$ up to a certain ideal $I_k$.

Consider the "geometric" action of $W_{aff}$ on $\htg^*$.
For any $\gamma \in \htg^*$ we can decompose
$\gamma = \lambda + c\delta + k\Lambda_0$, $\lambda \in
\gtc{h}^*$.
As $w \cdot \delta = \delta$ for any $w \in W_{aff}$ we can write
$$w\circ (\lambda + c\delta + k\Lambda_0) = w(\lambda,k) + (c+
d_w(\lambda, k))\delta + k\Lambda_0,$$
where  $w(\lambda,k) \in \gtc{h}^*$ and
$d_w(\lambda,k) = \left< d, w\circ(\lambda + k\Lambda_0) \right> \in
\bbb{Z}$.
Clearly, $w(\cdot, k)$ is a well--defined action of $W_{aff}$
on $\gtc{h}^*$.

Define $I_k$ as the linear subspace formed by $[\pi_\eta] -
(-1)^{l(w)} [\pi_{w(\eta,k)}]$ for $w \in W_{aff}$ and
dominant integral $\eta$, $w(\eta,k)$.

\begin{statement}{Proposition}

\rf $I_k$ is an ideal in $\ver{}$.

\rr The quotient $\ver{k}$ has a basis $\{ [\pi_\eta] \mid
\eta\in \Delta_k \}$.

\end{statement}

\begin{proof}

To prove (i) define the ${\cal{V}}(\gtc{g})$--module $\nmod$ with
basis $\{ [\eta]\}$, where $\eta$ is integral but not always dominant
weight. Let
$\pi$ be a $\gtc{g}$--module. Then we can write
$$[\pi] \cdot [\eta] = \sum_{\mu \in \gtc{h}^*} \dim (\pi_\mu)\cdot [\mu
+\eta],$$
where $\pi_\mu$ denotes the subspace of $\pi$ with weight $\mu$.

Clearly, it defines an action of ${\cal{V}}(\gtc{g})$ on $\nmod$.

The Weyl groups $W$ and $W_{aff}$ act on $\gtc{h}^*$ and, therefore, on
$\nmod$. Consider in $\nmod$ the subspace $J$ formed by elements $[\eta] -
(-1)^{l(w)}[w\circ \eta]$, $w \in W$ and the subspace $\widetilde{I}_k$
formed by elements $[\eta] - (-1)^{l(w)}[w(\eta,k)]$, $w \in W_{aff}$.

First of all, as the characters of $\gtc{g}$--modules are $W$--invariant,
$J$ is a submodule of $\nmod$. As linear parts of transformations
$w(\cdot,k)$, $w \in W_{aff}$ are elements of $W$, we have that
$\widetilde{I}_k$ is also a submodule of $\nmod$.

Then note that $\nmod / J$ is isomorphic to ${\cal{V}}(\gtc{g})$ as
${\cal{V}}(\gtc{g})$--modules. Clearly, they are isomorphic as
vector spaces, and one can
easily compare the actions of ${\cal{V}}(\gtc{g})$ using the Weyl formula
for character.

At last, we have $\widetilde{I}_k / J \cong I_k$.
So $I_k$ is a left ideal in ${\cal{V}}(\gtc{g})$. As
this algebra is commutative, we have (i).

The part (ii) is clear.

\end{proof}

Now we can formulate the rule for the dimensions of the spaces of
coinvariants.

\begin{statement}{Theorem}\label{mcoi}
Let ${\cal{P}} = (p_1, \dots, p_n)$.

\rf Let
$$ N_k = \bigoplus_{\eta \in \Delta_k} \pi_\eta.$$
We can write in $\ver{k}$ that
$$ [N_k]^n = \sum_{\eta \in \Delta_k} c_\eta [\pi_\eta].$$
Then $\dim V(k, \lambda, {\cal{P}}, \np) = c_\lambda$.

\rr Let
$$ NN_k = \bigoplus_{\eta \in \Delta_k} \pi_\eta \boxtimes \pi_\eta$$
be a $\gtc{g} \oplus \gtc{g}$--module. We can write in
$\ver{k}\otimes \ver{}$
$$ [NN_k]^n = \sum_{\eta \in \Delta_k;\, \mu \ \mbox{is dominant}}
c_\eta(\mu) \cdot [\pi_\eta] \otimes [\pi_\mu].$$
Then
$$V (k, \lambda, {\cal{P}}) = \bigoplus_{\mu} c_\lambda(\mu) \cdot \pi_\mu.$$

\end{statement}

Using the definition of the Verlinde algebra we can reformulate the theorem
for our purposes in the following way. If 
$$(N_k)^n = \bigoplus_{\eta} c(\eta) \cdot \pi_\eta$$
then
$$\dim V(k,\lambda, {\cal{P}}, \np) = \sum_\subw (-1)^{l(w)} c(w(\lambda,k)).$$

And if
$$(NN_k)^n = \bigoplus_{\eta, \mu} c(\eta, \mu) \cdot \pi_\eta \boxtimes
\pi_\mu$$
then
$$ V(k, \lambda, {\cal{P}}) = \sum_{\mu} \pi_\mu \cdot \sum_\subw
(-1)^{l(w)} c(w(\lambda,k), \mu).$$

\begin{proof}

To prove this theorem we need an interpretation of the Verlinde
algebra in terms of coinvariants.

Let as above
${\cal P}=\{p_1, \dots, p_n\}$ be a set of distinct points, let $l$ be
an integer. Let $\gtc{g}({\cal P}^l)=\gtc{g}\otimes \Complex
[t]\bigl(\prod (t-p_i)\bigr)^l$. The spaces $S_l=L_{k, \lambda}/\gtc{g}(
{\cal P}^l)$ form a projective system $S_0\mapsto S_1\mapsto
\dots$. Therefore, we can define the direct limit of dual spaces
$$
S=\lim_{\to}S_i^*.
$$
It has a natural structure of a $(\tg \oplus\dots \oplus \tg)$ ($n$
copies)--module, where $i$--th copy of $\tg$ is "settled" in
the point $p_i$.

\begin{statement}{Theorem}\label{refo}

\rf We have
\begin{equation}
S=\bigoplus_{\lambda_j \in \Delta_k}
C(\lambda_1, \dots, \lambda_n; \lambda) \cdot L_{\lambda_1,
k}\boxtimes \dots\boxtimes L_{\lambda_n, k}
\end{equation}

\rr Coefficients $C(\lambda_1, \dots, \lambda_n; \lambda)$ are "structural
constants" of the Verlinde algebra. Namely,

\begin{equation}
[\pi_{\lambda_1}]\cdot\dots\cdot [\pi_{\lambda_n}]=\sum_{\lambda\in \Delta_k}
C(\lambda_1, \dots, \lambda_n; \lambda) \cdot \pi_\lambda.
\end{equation}
\end{statement}

This statement can be treated as another definition of the Verlinde
algebra. For proof of (i) see \cite{BFM}, some methods how to
prove (ii) can be found in \cite{BFM}, \cite{F}.

The theorem \gref{appl}{mcoi} can be deduced from the theorem
\gref{appl}{refo}. Namely, one can show that the space
$V(k,\lambda,{\cal{P}}, \np)$
is dual to the space of highest weight vectors
in the decomposition of $S$ and that the space
$V(k,\lambda,{\cal{P}})$ is dual to the
$\gtc{g}\oplus\dots\oplus\gtc{g}$--module generated by
these hightest weight vectors.

\end{proof}

\subsection{Main conjecture.}

Now remark that $N_k$ is a cyclic $\gtc{g}$--module. Let $v_k$ be
the sum of highest vectors of irreducible components $\pi_\eta \subset
N_k$. Then $v_k$ is a cyclic vector of $N_k$.
Consider the filtered tensor product
${\cal{F}}_{\cal{Z}} (N_k, \dots, N_k)$ for general ${\cal{Z}}$ (see
section~\ref{init}). As the filtration is
$\gtc{g}$--equivariant, we can rewrite the
Hilbert polynomial of ${\cal{F}}_{\cal{Z}} (N_k, \dots,
N_k)$ in the following way:

$$ \ch_q(N_k^n) = \sum_\eta c_q(\eta) \cdot \ch \pi_\eta,$$

where $c_q(\eta)$ are polynomials in $q$.

\begin{statement}{Conjecture}\label{c1}
Let ${\cal{H}}_q(k,\lambda, {\cal{P}},\np)  = \sum q^i \dim \rm{gr}^i
V(k, \lambda, {\cal{P}}, \np)$. Then

$${\cal{H}}_q(k,\lambda, {\cal{P}},\np) = \sum_\subw (-1)^{l(w)}
q^{d_w(\lambda,k)} c_q(w(\lambda,k)).$$
Recall that $d_w(\lambda, k) = \left< d, w\circ (\lambda + k\Lambda_0)
\right>$.
\end{statement}

Now let us write the formula for $q$--character of this space.

The standard $\gtc{h}$--grading on the filtered tensor product doesn't
help
us to write the filtration, so consider the non--standard one.

Namely, let $\deg \pi_\eta = \eta$.
It means that the grading of each irreducible component coincides with
its highest weight.

This grading induces the grading on the tensor powers $(N_k)^n$.
One can easily show that the
filtration ${\cal{F}}_{\cal{Z}}(N_k, \dots , N_k)$,
structure of $\gtc{g}$--module and the grading on $(N_k)^n$ are
compatible.

So we can write
$$\ch_q^{\gtc{h}}(N_k^n) = \sum_\eta c_q^{\gtc{h}}(\eta) \cdot \ch
\pi_\eta.$$

\begin{statement}{Conjecture}\label{c2}
$$\ch_q V(k,\lambda, {\cal{P}},\np) = \sum_\subw (-1)^{l(w)}
q^{d_w(\lambda,k)} c_q^{\gtc{h}}(w(\lambda,k)).$$
\end{statement}

Now consider the space $V(k, \lambda, {\cal{P}})$. Similarly,
note that $NN_k$ is a cyclic $\gtc{g} \oplus \gtc{g}$--module.
Choose the sum of hightest vectors with respect to the action of
$\gtc{g} \oplus \gtc{g}$ as a cyclic vector and consider the filtered
tensor powers of $NN_k$. We can write
 $$\ch_q (NN_k^n) = \sum_{\eta,\mu} c_q(\eta,\mu) \cdot \ch
(\pi_\eta \boxtimes \pi_\mu).$$

\begin{statement}{Conjecture}\label{c3}
$$\ch_q V(k,\lambda, {\cal{P}}) = \sum_\mu \ch \, \pi_\mu \cdot
\sum_\subw (-1)^{l(w)}
q^{d_w(\lambda,k)} c_q (w(\lambda,k),\mu).$$
\end{statement}

\begin{Statement}{Remark}
These statements resemble the "parabolic" analog of the Weyl
formula for character. Namely, let $M_{k, \lambda}$ be the Weyl module.
By $q$ denote the variable responsible for the grading ''by $t$ ''.Then we
have
$$\ch L_{k, \lambda} = \sum_\subw (-1)^{l(w)} q^{d_w(\lambda,k)}
\ch M_{k, w(\lambda,k)}.$$
\end{Statement}

\noindent{\bf B.~Feigin}, Landau Institute of theor. physics,
Chernogolovka, Mosk. obl., Russia, e--mail {\bf
feigin@landau.ac.ru}

\medskip

\noindent{\bf S.~Loktev}, Independent University of Moscow, Bolshoj
Vlasievsky per., 7, Moscow, Russia, e--mail {\bf loktev@mccme.ru}

\end{document}